\documentclass[a4paper,12pt]{article}
\usepackage{amsthm,amsmath}
\usepackage{amssymb}
\usepackage[all]{xypic}
\usepackage{enumerate}
\usepackage{a4wide}
\usepackage{hyperref}
\usepackage{latexsym}
\usepackage{enumerate}

\newtheorem{theorem}{Theorem}[section]
\newtheorem{proposition}[theorem]{Proposition}
\newtheorem{corollary}[theorem]{Corollary}
\newtheorem{lemma}[theorem]{Lemma}
\newtheorem*{theorem*}{Theorem}
\newtheorem*{proposition*}{Proposition}
\newtheorem*{corollary*}{Corollary}
\newtheorem*{lemma*}{Lemma}
\theoremstyle{definition}
\newtheorem{definition}[theorem]{Definition}

\newtheorem{example}[theorem]{Example}
\newtheorem{remark}[theorem]{Remark}
\newtheorem*{remark*}{Remark}

\newtheorem*{definition*}{Definition}

\newcommand{\cat}[1]{\mathcal{#1}}
\newcommand{\coring}[1]{\mathfrak{#1}}
\newcommand{\tensor}[1]{\otimes_{#1}}
\newcommand{\tensfun}[1]{\underset{{#1}}{\otimes}}

\newcommand{\rcomod}[1]{\mathcal{M}^{#1}}
\newcommand{\rmod}[1]{\mathcal{M}_{#1}}
\newcommand{\lmod}[1]{{}_{#1}\mathcal{M}}

\newcommand{\cotensor}[1]{\square_{#1}}

\renewcommand{\hom}[3]{\mathrm{Hom}_{#1}(#2,#3)}

\newcommand{\End}[2]{\mathrm{End}_{#1}(#2)}

\newcommand{\rend}[2]{\mathrm{End}({#2}_{#1})}
\newcommand{\lend}[2]{\mathrm{End}({}_{#1}#2)}

\newcommand{\rcomatrix}[2]{#2^* \tensor{#1} #2}

\newcommand{\add}[1]{\mathrm{add}(#1)}

\newcommand{\Rcn}[1]{\mathfrak{R}(#1)}
\newcommand{\ffun}[1]{{}^{\omega}#1}

\begin{document}
\title{Infinite Comatrix Corings}
\author{L. El Kaoutit {\normalsize and} J. G\'omez-Torrecillas \\
\normalsize Departamento de \'{A}lgebra \\ \normalsize Universidad
de Granada\\ \normalsize E18071 Granada, Spain \\ \normalsize
e-mail: \textsf{kaoutit@fedro.ugr.es} \\ \normalsize e-mail:
\textsf{torrecil@ugr.es} }

\date{\empty}

\maketitle

\begin{abstract}
We characterize the corings whose category of comodules has a generating set of small
projective comodules in terms of the (non commutative) descent theory. In order to extricate
the structure of these corings, we give a generalization of the notions of comatrix coring and
Galois comodule which avoid finiteness conditions. A sufficient condition for a coring to be
isomorphic to an infinite comatrix coring is found. We deduce in particular that any coalgebra
over a field and the coring associated to a group-graded ring are isomorphic to adequate
infinite comatrix corings. We also characterize when the free module canonically associated to
a (not necessarily finite) set of group like elements is Galois.
\end{abstract}

\section*{Introduction}

Among the different aspects in the recent developments of the theory of corings, one of the
most intensively studied is the notion of a Galois coring, and its relationships with
(noncommutative) descent theory for ring extensions and Morita-type equivalence theorems. A
coring $\coring{C}$ over a ring $A$ is said to be Galois \cite{Brzezinski:2002} whenever the
canonical map $\mathsf{can}: A \tensor{B} A \rightarrow \coring{C}$ which sends $a \tensor{B}
a'$ onto $aga'$ is an isomorphism, where $g \in \coring{C}$ is a grouplike element and $B = \{
a \in A \; | \; ag = ga \}$. One of the origins of this notion is the concept of a
noncommutative Hopf-Galois extension \cite{Kreimer/Takeuchi:1981,
Ulbrich:1982,Doi/Takeuchi:1989, Schneider:1990}, which can be ultimately traced back to the
theory of principal homogeneous spaces for actions of affine groups over affine schemes (see,
e.g. \cite[18.3]{Waterhouse:1979}), and the characterization of strongly graded rings
\cite{Dade:1980}. One of the fundamental facts is that a faithfully flat Hopf-Galois extension
$B \subseteq A$, for $H$ a Hopf algebra, encodes a canonical equivalence of categories between
a category of Hopf modules and the category of right modules over $B$ (see
\cite{Schneider:1990}). This has been a model for research in more general frameworks (see
\cite{Caenepeel:2003unp} for a helpful survey), including coalgebra-Galois extensions for
entwining structures \cite{Brzezinski/Hajac:1999}, and Galois corings with a grouplike
\cite{Brzezinski:2002}. In \cite{ElKaoutit/Gomez:2003} it is shown that, in order to formulate
a meaningful notion of Galois coring, the grouplike can be replaced by a right comodule (a
Galois comodule in the terminology of \cite{Brzezinski/Wisbauer:2003,Brzezinski:unp2003}),
finitely generated and projective as a module over the ground ring $A$, whenever the role of
Sweedler's canonical coring is played by the associated comatrix coring. This new viewpoint
led us to imagine that there were a relation between the notion of Galois coring and the fact
that a coalgebra over a field can be reconstructed from its finite dimensional comodules
\cite{Joyal/Street:1990}. This paper contains the mathematical results of our investigations
on this idea. A generalization of the notion of comatrix coring (and of Galois comodule) will
be introduced with this purpose. The role of Galois comodules in Non Commutative Geometry, as
non commutative principal bundles, is revealed in \cite{Brzezinski:unp2003}.

We will define our generalized comatrix corings in three ways, being the interplay between
these constructions fundamental in our exposition. One of them is inspired in the coalgebra
defined in \cite[Section 4]{Joyal/Street:1990}. We wish to thank to Antonio M. Cegarra for
helping us to understand this construction.

\section{A specialized introduction: Galois corings with
several grouplike elements}

In this section, we will give our statements without proofs. They
are consequences of the results proved in the rest of the paper.
The basic notations and notions will explained later, as well.

Let $G$ be a set of grouplike elements in a coring $\coring{C}$
over a $K$--algebra $A$ ($K$ is a commutative ring). For each $g
\in G$, we will denote by $[g]A$ the right $\coring{C}$--comodule
structure on $A$ with coaction $A \rightarrow A \tensor{A}
\coring{C}$ sending $a \in A$ onto $1 \tensor{A} ga$. The notation
$A[g]$ stands for the left $\coring{C}$--comodule structure on $A$
associated to $g$. For $g, h \in G$, the $K$--module
$\hom{\coring{C}}{[g]A}{[h]A}$ is identified, via the map $f
\mapsto f(1)$, with $A_{g,h} = \{ b \in A ~|~ hb = bg \}$.
Clearly, $A_{g,g}$ is a (unital) subring of $A$, and $A_{g,h}$
becomes an $A_{h,h}-A_{g,g}$--bimodule. Consider the ring (not
necessarily with unit) $S = M_{G \times G}^f$ consisting of all $G
\times G$--matrices with coefficients in $A$ and with finitely
many non-zero entries. The unit $1 \in A$ at the position $(g,g)$
and zero elsewhere gives an idempotent matrix $1_{g,g} \in S$, and
the set of all $1_{g,g}$'s gives a complete set of pairwise
ortogonal idempotents for $S$. The external direct sum $R =
\bigoplus_{g,h \in G}A_{g,h}$ may be then considered as a subring
of $S$ that contains the idempotents $1_{g,g}$.

On the other hand, for every $g \in G$ we can consider the
Sweedler's canonical $A$--coring $A \tensor{A_{g,g}} A$
\cite{Sweedler:1975}. Let $\mathfrak{J}$ be the $K$--submodule of
the coproduct $\bigoplus_{g \in G} A \tensor{A_{g,g}} A$ generated
by all elements of the form $a \tensor{A_{h,h}} ta' - at
\tensor{A_{g,g}} a'$, where $g, h \in G$, $a, a' \in A$, and $t
\in A_{g,h}$. It turns out that $\mathfrak{J}$ is a coideal of the
coproduct $A$--coring $\bigoplus_{g \in G} A \tensor{A_{g,g}} A$,
and we have a new $A$--coring
\begin{equation*}
\mathfrak{R}(G) = \frac{\bigoplus_{g \in G} A \tensor{A_{g,g}}
A}{\mathfrak{J}} = \frac{\bigoplus_{g \in G} A \tensor{A_{g,g}}
A}{\langle a \tensor{A_{h,h}} ta' - at \tensor{A_{g,g}} a'; \;
a,a' \in A, \; t \in A_{g,h}, \; g,h \in G \rangle}
\end{equation*}

Moreover, the map $\mathsf{can} : \mathfrak{R}(G) \rightarrow \coring{C}$ which sends $a
\tensor{A_{g,g}} a' + \mathfrak{J}$ onto $aga'$ is a homomorphism of $A$--corings. Obviously,
this canonical map is a generalization of the given in \cite{Brzezinski:2002} for a single
group-like (i.e., $G = \{ g \}$). The comodule $\Sigma = \bigoplus_{g \in G}[g]A$ is not
finitely generated unless $G$ is finite. Nevertheless, following \cite{Brzezinski:2002} for
the case $G = \{g\}$, it makes sense to say that $\Sigma$ is a Galois coring (with respect to
$\Sigma$) whenever $\mathsf{can}$ is an isomorphism. Of course, there are examples of
situations where $\mathsf{can}$ is bijective but $\Sigma$ is not finitely generated. The key
is the following result.

\begin{theorem}\label{galoisG}
If $\bigoplus_{g \in G}[g]A$ is a generator for
$\rcomod{\coring{C}}$ then $\mathsf{can}: \mathfrak{R}(G)
\rightarrow \coring{C}$ is an isomorphism of $A$--corings.
\end{theorem}

A relevant example of Galois coring in this new general framework is the following.

\begin{example}\label{graduados}
Let $A$ be a $G$--graded ring, where $G$ is any group. Endow the free left $A$--module $AG$
with basis $G$ with the right $A$--module structure given by $ga = agh$, for $a \in A$
homogeneous of degree $h \in G$. Then $AG$ becomes an $A$--bimodule. Consider the $A$--coring
structure on $AG$ given by the comultiplication defined by $\Delta (ag) = ag \tensor{A} g$ and
counit by $\epsilon(ag) = a$, for every $a \in A$, $g \in G$. The category of right
$AG$--comodules is then isomorphic with the category $gr-A$ of all $G$--graded right
$A$--modules. This can be proved either by direct computations or by using that $AG$ is the
coring built, according to \cite[Proposition 2.2]{Brzezinski:2002} and \cite[Example
3.1]{Brzezinski:1999}, from the entwining structure given on $AG \cong A \tensor{K} KG$ by the
$KG$--comodule algebra map $a \mapsto \sum_{g \in G} a_g \tensor{} g$, and the fact that, in
this case, the category of Hopf modules $\mathcal{M}_A^{KG}$ is isomorphic $gr-A$
\cite{Koppinen:1994}. Clearly, $G$ is a set of group like elements, and $\{ [g]A : g \in G \}$
is nothing but the set of all shifts of the graded module $A$. It is known that it is a
generating set of small projectives for $gr-A$. By Theorem \ref{galoisG}, the canonical map is
an isomorphism. More generally, for each subgroup $H$ of $G$ such that $\bigoplus_{h \in
H}[h]A$ is a (projective) generator for $gr-A$, the canonical map $\mathsf{can} :
\mathfrak{R}(H) \rightarrow AG$ is an isomorphism of $A$--corings.
\end{example}

\begin{theorem}\label{muchosgrouplike}
The following statements are equivalent for an $A$--coring
$\coring{C}$ with a set of group-like elements $G$.
\begin{enumerate}[(i)]
\item ${}_A\coring{C}$ is flat, $can:\mathfrak{R}(G) \rightarrow
\coring{C}$ is an isomorphism, and $S$ is faithfully flat as a left $R$--module; \item
${}_A\coring{C}$ is flat, $can : \mathfrak{R}(G) \rightarrow \coring{C}$ is an isomorphism and
the category $\rcomod{\mathfrak{R}(G)}$ of right $\mathfrak{R}(G)$--comodules is equivalent,
in a canonical way, to the category $\rmod{R}$ of all unital right $R$--modules; \item
${}_A\coring{C}$ is flat and $\bigoplus_{g \in G} [g]A$ projective generator for
$\rcomod{\coring{C}}$;
\item ${}_A\coring{C}$ is flat and the category $\rcomod{\coring{C}}$ is equivalent, in a
canonical way, to the category $\rmod{R}$ of all unital right $R$--modules.
\end{enumerate}
\end{theorem}

We will in fact prove a more general result (Theorem
\ref{descent}), characterizing those corings whose category of
comodules has a generating set of small projectives. With this
purpose, we extend the techniques developed in
\cite{ElKaoutit/Gomez:2003} (in particular the notion of a
comatrix coring) to a larger framework. We thus will obtain an
alternative description of the coring $\mathfrak{R}(G)$ as
follows: write ${\Sigma}^{\dagger} = \bigoplus_{g \in G} A[g]$,
$\Sigma = \bigoplus_{g \in G}[g]A$. Let us understand $\Sigma$
(resp. ${\Sigma}^{\dagger}$) as the free right (resp. left)
$A$--module with basis $\{ [g] : g \in G \}$. Consider the
$A$--coring ${\Sigma}^{\dagger} \tensor{R} \Sigma$ with
comultiplication given by
\[
\Delta(\sum_{g \in G} a_g[g] \tensor{R} \sum_{h \in G} [h]a'_{h})
= \sum_{g \in G} a_g[g] \tensor{R} [g]1 \tensor{A} 1[g] \tensor{R}
[g]a'_g,
\]
and counit
\[
\epsilon(\sum_{g \in G} a_g[g] \tensor{R} \sum_{h \in G}
[h]a'_{h}) = \sum_{g \in G}a_ga'_g
\]
Then there exists a canonical isomorphism of $A$--corings $\mathfrak{R}(G) \cong
{\Sigma}^{\dagger} \tensor{R} \Sigma$, and this last is a sort of generalized comatrix coring.
In fact, this setting allows to extend the methods from \cite{ElKaoutit/Gomez:2003} to a more
general context. As a consequence, we will have that the corings characterized in Theorem
\ref{muchosgrouplike} are those for which ${}_R\Sigma$ is faithfully flat.

\begin{example}\label{graduados1}
Continuing with Example \ref{graduados} we get that the category $gr-A$ is always equivalent
to the category of unital right modules over the ring $R$. More generally, given a subgroup
$H$ of $G$, Theorem \ref{muchosgrouplike} gives several conditions which characterize when
$gr-A$ is equivalent, in a canonical way, to the category of unital right $R$--modules where
$R$ is the ring of $H \times H$--matrices with finitely many nonzero entries whose
coefficients in its $(g,h)$--entry are the homogeneous elements of $A$ of degree $gh^{-1}$,
where $g, h \in H$.
\end{example}

\section{Basic notions}

We use the following conventions. We work over a fixed commutative
ring $K$, and all our additive categories are assumed to be
$K$--linear. For instance, all rings in this paper are (not
necessarily unitary) $K$--algebras, and all bimodules are assumed
to centralize the elements of $K$. Associated to every object $C$
of an additive category $\mathcal{C}$ we have its endomorphism
ring $\mathrm{End}_{\mathcal{C}}(C)$, whose multiplication is
given by the composition of the category. As usual, some special
conventions will be understood for the case of endomorphism rings
of modules. Thus, if $M_R$ is a (unital) right module over a ring
$R$, then its endomorphism ring in the category $\rmod{R}$ of all
unital right $R$--modules will be denoted by $\rend{R}{M}$, while
if ${}_RN$ is a left $R$--module, then its endomorphism ring,
denoted by $\lend{R}{N}$, is, by definition, the opposite of the
endomorphism ring of $N$ in the category $\lmod{R}$ of all unital
left modules over $R$. We will make use of rings $R$ which need
not to have a unit, although they will always contain a set of
pairwise ortogonal idempotents $\{ e_i \}$ which is complete in
the sense that $R = \bigoplus_{i} Re_i = \bigoplus_{i}e_iR$. To be
unital for a right $R$--module $M$ means that $M = MR$. The tensor
product over $R$ is denoted by $\tensor{R}$. We shall sometimes
replace $\tensor{K}$ by $\tensor{}$.

The notation $A$ is reserved for a $K$--algebra with unit. We
recall from \cite{Sweedler:1975} that an $A$--\emph{coring} is a
trhee-tuple ($\coring{C}, \Delta,\epsilon$) consisting of an
$A$-bimodule $\coring{C}$ and two $A$--bimodule maps
\[
\Delta : \coring{C} \longrightarrow \coring{C} \tensor{A} \coring{C}, \qquad \epsilon :
\coring{C} \longrightarrow A
\]
such that $(\coring{C} \tensor{A} \Delta) \circ \Delta = (\Delta
\tensor{A} \coring{C}) \circ \Delta$ and $(\epsilon \tensor{A}
\coring{C}) \circ \Delta = ( \coring{C} \tensor{A} \epsilon) \circ
\Delta = 1_\coring{C}$. A \emph{right}
$\coring{C}$--\emph{comodule} is a pair $(M,\rho_M)$ consisting of
a right $A$--module $M$ and an $A$--linear map $\rho_M: M
\rightarrow M \tensor{A} \coring{C}$ satisfying $(M \tensor{A}
\Delta) \circ \rho_M = (\rho_M \tensor{A} \coring{C}) \circ
\rho_M$, $(M \tensor{A} \epsilon) \circ \rho_M = 1_M$. A
\emph{morphism} of right $\coring{C}$--comodules $(M,\rho_M)$ and
$(N,\rho_N)$ is a right $A$--linear map $f: M \rightarrow N$ such
that $(f \tensor{A} \coring{C}) \circ \rho_M = \rho_N \circ f$;
the $K$--module of all such morphisms will be denoted by
$\hom{\coring{C}}{M}{N}$. The right $\coring{C}$--comodules
together with their morphisms form the additive category
$\rcomod{\coring{C}}$. Coproducts and cokernels in
$\rcomod{\coring{C}}$ do exist and can be already computed in
$\rmod{A}$. Therefore, $\rcomod{\coring{C}}$ has arbitrary
inductive limits. If ${}_A\coring{C}$ is flat, then
$\rcomod{\coring{C}}$ is an abelian category. The converse is not
true, as \cite[Example 1.1]{ElKaoutit/Gomez/Lobillo:2001unp}
shows.

 Let $\rho_M : M
\rightarrow M \tensor{A} \coring{C}$ be a comodule structure over an $A'-A$--bimodule $M$, and
assume that $\rho_M$ is $A'$--linear ($A'$ denotes a second unitary $K$-algebra). For any
right $A'$--module $X$, the right $A$--linear map $X \tensor{A'} \rho_{M} : X \tensor{A'} M
\rightarrow X \tensor{A'} M \tensor{A} \coring{C}$ makes $X \tensor{A'} M$ a right
$\coring{C}$--comodule. This leads to an additive functor $- \tensor{A'} M : \rmod{A'}
\rightarrow \rcomod{\coring{C}}$. The classical adjointness isomorphism $\hom{A}{Y
\tensor{A'}M}{X} \cong \hom{A'}{Y}{\hom{A}{M}{X}}$ induces, by restriction, a natural
isomorphism $\hom{\coring{C}}{Y \tensor{A'}M}{X} \cong \hom{A'}{Y}{\hom{\coring{C}}{M}{X}}$,
for $Y \in \rmod{A'}, X \in \rcomod{\coring{C}}$. Therefore, $\hom{\coring{C}}{M}{-} :
\rcomod{\coring{C}} \rightarrow \rmod{A'}$ is right adjoint to $- \tensor{A'} M : \rmod{A'}
\rightarrow \rcomod{\coring{C}}$. On the other hand, the functor $- \tensor{A} \coring{C}$ is
right adjoint to the forgetful functor $U : \rmod{A} \rightarrow \rcomod{\coring{C}}$ (see
\cite[Proposition 3.1]{Guzman:1989}, \cite[Lemma 3.1]{Brzezinski:2002}).

Now assume that the $A'-A$--bimodule $M$ is also a left comodule over an $A'$--coring
$\coring{C}'$ with structure map $\lambda_M : M \rightarrow \coring{C}' \tensor{A} M$. It is
clear that $\rho_M : M \rightarrow M \tensor{A} \coring{C}$ is a morphism of left
$\coring{C}'$--comodules if and only if $\lambda_M : M \rightarrow \coring{C}' \tensor{A'} M$
is a morphism of right $\coring{C}$--comodules. In this case, we say that $M$ is a
$\coring{C}'-\coring{C}$--bicomodule.


A grouplike is an element $g \in \coring{C}$ such that $\Delta(g) = g \tensor{A} g$ and
$\epsilon(g) = 1$. Every grouplike defines a right $\coring{C}$--comodule structure $\rho : A
\rightarrow A \tensor{A} \coring{C} \cong \coring{C}$ given by $\rho(a) = ga$. A left comodule
structure is similarly obtained. This process can be reversed: each right
$\coring{C}$--comodule structure $\rho : A \rightarrow \coring{C}$ gives a grouplike $g =
\rho(1)$ (see \cite[Lemma 5.1]{Brzezinski:2002}).

A friendly source of information on corings and their comodules is the monograph
\cite{Brzezinski/Wisbauer:2003}.

\section{Reconstruction and infinite comatrix corings}\label{inf-matrix}

A ring extension $B \subseteq A$ for rings with unit can be understood as a faithful functor
between two additive categories with one object. With this idea in mind, and with an eye on
\cite[Section 4]{Joyal/Street:1990} the construction of the canonical $A$--coring $A
\tensor{B} A$ from \cite{Sweedler:1975} is generalized in this section.

 Let $K$ denote a commutative ring.
All categories and functors will be assumed to be $K$--linear. Let $\add{A_A}$ denote the
category of all finitely generated projective right modules over an associative $K$--algebra
with unit $A$. Let $\omega : \cat{A} \rightarrow \add{A_A}$ be a functor, where $\cat{A}$ is a
small category. The image of an object $P$ of $\cat{A}$ under $\omega$ will be denoted by
$\ffun{P}$, or even by $P$ itself, when no confusion may be expected. For each $P \in
\cat{A}$, there is a canonical homomorphism of $K$--algebras from $T_P = \End{\cat{A}}{P}$ to
$S_P = \rend{A}{\ffun{P}}$, which sends $t$ onto $\omega(t)$. The module $\ffun{P}$ becomes
then an $T_P-A$--bimodule. If $P$, $Q$ are objects of $\cat{A}$, then the
 elements of the $T_Q-T_P$--bimodule $T_{PQ} =
 \hom{\cat{A}}{P}{Q}$ act canonically on $\ffun{P}$. This action
 can be thought as the $T_Q-A$--bimodule map
\[
\xymatrix@R=0pt@C=50pt{ T_{PQ} \tensor{T_P} \ffun{P} \ar@{->}[r] &
\ffun{Q}
\\ t \tensor{T_P} p \ar@{|->}[r] & tp:=
\omega(t)(p)}
\]
The right dual $A$--modules $\ffun{P}^*= \hom{A}{\ffun{P}}{A}$ are
then in a natural way  $A-T_P$--bimodules; the corresponding dual
pairings are given by
\[
\xymatrix@R=0pt@C=50pt{ \ffun{Q}^* \tensor{T_Q} T_{PQ} \ar@{->}[r]
& \ffun{P}^* \\ \varphi \tensor{T_Q} t \ar@{|->}[r] & \varphi t :=
\varphi \circ \omega(t) }
\]
From now on, we will denote $\ffun{P}$ by $P$. We can associate to every object $P$ of
$\cat{A}$ its \emph{comatrix $A$--coring} $\rcomatrix{T_P}{P}$, defined as follows
\cite{ElKaoutit/Gomez:2003}. Let $\{ e_{\alpha_P}^*, e_{\alpha_P} \}$ be a finite dual basis
for the module $P_A$, and define the comultiplication as $\Delta (\varphi \tensor{T_P} p) =
\sum_{\alpha_P}\varphi \tensor{T_P} e_{\alpha_P} \tensor{A} e_{\alpha_P}^* \tensor{T_P} p$.
The counit is given by $\epsilon (\varphi \tensor{T_P} p) = \varphi (p)$. We can then consider
the coproduct of $A$--corings
\[
\mathfrak{P}(\cat{A}) = \bigoplus_{P \in \cat{A}}
\rcomatrix{T_P}{P}
\]
Every $P \in \cat{A}$ is canonically a right $\rcomatrix{T_P}{P}$--comodule
\cite{ElKaoutit/Gomez:2003} and, hence, a $\mathfrak{P}(\cat{A})$--comodule, with structure
map $\rho_P : P \rightarrow P \tensor{A} \mathfrak{P}(\cat{A})$ defined as $\rho_P (p) =
\sum_{\alpha_P} p \tensor{A} e_{\alpha_P}^* \tensor{T_P} e_{\alpha_P}$. The assignment $P
\mapsto (P, \rho_P)$ is not, at least in the obvious way, a functor from $\cat{A}$ to
$\rcomod{\mathfrak{P}(\cat{A})}$. In order to remedy this, we will factor out
$\mathfrak{P}(\cat{A})$ by a coideal.

\begin{lemma}\label{compatible}
Let $t \in T_{PQ}$. Then
\[
\sum_{\alpha_{Q}} e_{\alpha_{Q}} \tensor{A} e_{\alpha_{Q}}^* t =
\sum_{\alpha_{P}} t e_{\alpha_{P}} \tensor{A} e_{\alpha_{P}}^*.
\]
\end{lemma}
\begin{proof}
Follows easily from the dual basis criterion.
\end{proof}

\begin{lemma}\label{Jota}
 The $K$--submodule $\mathfrak{J}$ of $\mathfrak{P}(\cat{A})$ generated by the set
\[
\{ \varphi\tensor{T_Q} tp - \varphi t \tensor{T_P} p \; : \;
\varphi \in Q^*, p \in P, t \in T_{PQ}, P,Q \in \cat{A} \}
\]
is a coideal of $\mathfrak{P}(\cat{A})$.
\end{lemma}
\begin{proof}It is easily shown that, in fact, $\mathfrak{J}$ is
an $A$--subbimodule of $\mathfrak{P}(\cat{A})$. Let us now prove
that $\epsilon(\mathfrak{J}) = 0$. In fact,
\[
\begin{array}{ll}
\epsilon(\varphi \tensor{T_{Q}} tp - \varphi
 t \tensor{T_{P}} p) & = \varphi(tp)-
\varphi  t ( p) \\ & = \varphi t(p) - \varphi t ( p) =0.
\end{array}
\]
Now, by using Lemma \ref{compatible}, we have

\begin{multline}
\Delta(\varphi \tensor{T_Q} tp -
\varphi  t \tensor{T_P} p) \\
= \sum_{\alpha_Q} \varphi \tensor{T_Q} e_{\alpha_Q} \tensor{A}
e_{\alpha_Q}^* \tensor{T_Q} tp - \sum_{\alpha_P} \varphi t
\tensor{T_P} e_{\alpha_P}
\tensor{A} e_{\alpha_P}^* \tensor{T_P} p \\
= \sum_{\alpha_Q} \varphi \tensor{T_Q} e_{\alpha_Q} \tensor{A}
e_{\alpha_Q}^* \tensor{T_Q} tp- \sum_{\alpha_Q} \varphi
\tensor{T_Q} e_{\alpha_Q}
\tensor{A} e_{\alpha_Q}^*  t \tensor{T_P} p  \\
+ \sum_{\alpha_P} \varphi \tensor{T_Q} te_{\alpha_P} \tensor{A} e_{\alpha_P}^* \tensor{T_P} p
- \sum_{\alpha_P} \varphi t \tensor{T_P} e_{\alpha_P} \tensor{A} e_{\alpha_P}^* \tensor{T_P} p
\\ = \sum_{\alpha_Q} \varphi \tensor{T_Q} e_{\alpha_Q} \tensor{A}
(e_{\alpha_Q}^* \tensor{T_Q} tp-
 e_{\alpha_Q}^*
t \tensor{T_P} p )\\
+ \sum_{\alpha_P} (\varphi \tensor{T_Q} te_{\alpha_P} - \varphi t
\tensor{T_P} e_{\alpha_P}) \tensor{A} e_{\alpha_P}^* \tensor{T_P}
p
\end{multline}
This proves that $\Delta(\mathfrak{J}) \subseteq Ker (\pi
\tensor{A} \pi)$, where $\pi : \mathfrak{P}(\cat{A}) \rightarrow
\mathfrak{P}(\cat{A})/\mathfrak{J}$ is the canonical projection.
Therefore, $\mathfrak{J}$ is a coideal.
\end{proof}

\begin{proposition}
Let $\mathfrak{P}(\cat{A}) = \bigoplus_{P \in \cat{A}}
\rcomatrix{T_P}{P}$ and define the factor $A$--coring
$\Rcn{\cat{A}} = \mathfrak{P}(\cat{A})/\mathfrak{J}$. There is a
functor $\Rcn{\omega_A} : \cat{A} \rightarrow
\rcomod{\Rcn{\cat{A}}}$ making the diagram
\begin{equation}\label{diag2bis}
\xymatrix@R=5pt@C=50pt{\cat{A} \ar@{->}^{\omega_A}[r]
\ar@{.>}_{\Rcn{\omega_A}}[dd] & \add{A_A} \ar[dd]\\
& \\
 \rcomod{\Rcn{\cat{A}}} \ar@{->}_{U_A}[r] & \rmod{A}}
\end{equation}
commutative. This functor assigns to every $P \in \cat{A}$ the
right $\Rcn{\cat{A}}$--comodule induced by its canonical
$\rcomatrix{T_P}{P}$--coaction and the canonical homomorphism of
$A$--corings $\pi : \mathfrak{P}(\cat{A}) \rightarrow
\mathfrak{P}(\cat{A})/\mathfrak{J}$.
\end{proposition}
\begin{proof}
So, the right $\Rcn{\cat{A}}$--comodule structure for $P$ is given
by
\[
\xymatrix{P \ar^-{\rho_{P}}[rr] & & P \tensor{A}
\rcomatrix{T_P}{P} \ar^-{P \tensor{A} \pi}[rr] & & P
\tensor{A}\Rcn{\cat{A}}\\
p \ar@{|->}[rrrr]& & & & \sum_{\alpha_P}e_{\alpha_P} \tensor{A}
(e_{\alpha_P}^* \tensor{T_P} p + \mathfrak{J})}
\]
Given $\lambda \in \hom{\cat{A}}{P}{Q}$, a straightforward
computation shows, with the help of Lemma \ref{compatible}, that
$\omega(\lambda)$ is $\Rcn{\cat{A}}$--colinear.
\end{proof}

We shall now give an alternative description of the $A$--coring $\mathfrak{R}(\cat{A})$.
Assume that $\cat{A}$ is a subcategory of an additive category $\cat{C}$, and that there
exists the coproduct $\Sigma = \bigoplus_{P \in \cat{A}} P$ in the category $\cat{C}$. Assume
further that the functor $\omega : \cat{A} \rightarrow \add{A_A}$ is the restriction of a
functor $U : \cat{C} \rightarrow \rmod{A}$ which commutes with the coproduct $\bigoplus_{P \in
\cat{A}}P$. This is not actually a restriction: such a category $\cat{C}$ can be constructed
by introducing formally a new object $\Sigma$, and enlarging the set of morphisms by defining
the new hom-sets $\hom{\cat{C}}{P}{Q} = \hom{\cat{A}}{P}{Q}$, $\hom{\cat{C}}{P}{\Sigma} =
\bigoplus_{X \in \cat{A}}\hom{\cat{C}}{P}{X}$, $\hom{\cat{C}}{\Sigma}{Q} = \prod_{X \in
\cat{A}} \hom{\cat{C}}{X}{Q}$ and, finally, $\hom{\cat{C}}{\Sigma}{\Sigma} = \prod_{X \in
\cat{A}}\hom{\cat{C}}{X}{\Sigma}$, for $P , Q \in \cat{A}$. The definition of the functor $U :
\cat{C} \rightarrow \rmod{A}$ is then clear.

Consider the endomorphism ring $T = \End{\cat{C}}{\Sigma}$. We have then a canonical structure
of a $T-A$--bimodule on $U(\Sigma)$, which induces an $A-T$--bimodule structure on
$U(\Sigma)^* = \hom{A}{\Sigma}{A}$. We think there will be no problems by using the notation
$\Sigma$ instead of $U(\Sigma)$. For each object $P$ in $\cat{A}$, let $\pi_P : \Sigma
\rightarrow P$ (resp. $\iota_P : P \rightarrow \Sigma$) be the canonical projection (resp.
injection).

\begin{proposition}\label{presentacion}
There is a surjective homomorphism of $A$--bimodules
\[
\Gamma : \bigoplus_{P \in \cat{A}} \rcomatrix{T_P}{P} \rightarrow \rcomatrix{T}{\Sigma}
\]
whose restriction $\Gamma_P$ to each $\rcomatrix{T_P}{P}$ is given by $\Gamma_P (\varphi
\tensor{T_P} p )= \varphi \pi_P \tensor{T} \iota_P (p)$. The kernel of $\Gamma$ is the coideal
$\mathfrak{J}$ and, hence, $\rcomatrix{T}{\Sigma}$ can be endowed with a structure of
$A$--coring such that $\Gamma$ is an homomorphism of $A$--corings. This homomorphism $\Gamma$
induces an isomorphism of corings $\mathfrak{R}(\cat{A}) \cong \rcomatrix{T}{\Sigma}$.
\end{proposition}
\begin{proof}
Let us first check that each $\Gamma_P$ is well defined. For
$\varphi \in P^*$, $p \in P$ and $t \in T_P$ we have
\[
\begin{array}{lclcl}
 \Gamma_P(\varphi t \tensor{T_P} p) & = & (\varphi t)\pi_P
 \tensor{T} \iota_P(p) & & \\
 & = & \varphi \pi_P \iota_P t \pi_P \tensor{T} p & &  \\
 & = & \varphi \pi_P \tensor{T} \iota_P t \pi_P \iota_P (p) & & \\
 & = & \varphi \pi_P \tensor{T} \iota_P t (p) & = &
 \Gamma_P(\varphi \tensor{T_P} t p)
\end{array}
\]
To prove that $\Gamma$ is surjective, observe that if $\varphi \tensor{T} p \in
\rcomatrix{T}{\Sigma}$ then there is a finite set $\cat{F}$ of objects of $\cat{A}$ such that
$p = (\sum_{P \in \cat{F}}\iota_P \pi_P)(p)$. Therefore
\begin{equation}\label{sobreyectivo}
\begin{array}{lclcl}
\varphi \tensor{T} p & = & \varphi \tensor{T} \sum_{P \in \cat{F}}\varphi \tensor{T} \iota_P \pi_P (p) & & \\
& = & \sum_{P \in \cat{F}} \varphi \iota_P \pi_P \tensor{T}
\iota_P \pi_P (p) & = & \Gamma(\sum_{P \in \cat{F}} \varphi
\iota_P \tensor{T_P} \pi_P (p))
\end{array}
\end{equation}
Next, let us check that $\mathfrak{J} \subseteq Ker \Gamma$. Given a generator
$\varphi\tensor{T_Q} tp - \varphi t \tensor{T_P} p$ of $\mathfrak{J}$, where $\varphi \in Q^*,
p \in P, t \in T_{PQ}$, and $P,Q \in \cat{A}$, apply $\Gamma$ to obtain
\[
\begin{array}{lclcl}
 \Gamma(\varphi \tensor{T_Q} tp - \varphi t \tensor{T_P} p) & =  &
 \varphi \pi_Q \tensor{T} \iota_Q(tp) - \varphi t \pi_P \tensor{T}
 \iota_P(p) & & \\
 & = & \varphi \pi_Q \iota_Q(tp) - \varphi \pi_Q \iota_Q t \pi_P
 \tensor{T} \iota_P (p) & & \\
 & = & \varphi \pi_Q \tensor{T} \iota_Q (tp) - \varphi \pi_Q
 \tensor{T} \iota_Q t \pi_P \iota_P (p) & & \\
 & = & \varphi \pi_Q \tensor{T} \iota_Q(t p) - \varphi \pi_Q
 \tensor{T} \iota_Q(t p) & = & 0
\end{array}
\]
Let us finally check that $Ker \Gamma \subseteq \mathfrak{J}$. An arbitrary element of
$\bigoplus_{P \in \cat{A}} \rcomatrix{T_P}{P}$ is a finite sum $x = \sum_{a = 1}^n \sum_{P \in
\mathcal{A}} \varphi_{a,P} \tensor{T_P} p_{a,P}$, for some $\varphi_{a,P} \in P^*$ and
$p_{a,P} \in P$ with almost all $p_{a,P} = 0$. Such an element becomes to $Ker \Gamma$ if and
only if $\sum_{a,P} \varphi_{a,P}\pi_{P} \tensor{T} \iota_P(p_{a,P}) = 0$. Since the elements
of the form $\iota_P(p)$ generate $\Sigma$, there exist \cite[Proposition
I.8.8]{Stenstrom:1975} a finite set $\{ \nu_k \}_{k \in I} \subseteq \Sigma^*$ and a set $\{
g_{a,P,k} \} \subseteq T$ such that
\begin{enumerate}
\item $g_{a,P,k} = 0$ for almost all $(a,P,k)$,
\item $\sum_{a,P}g_{a,P,k} \iota_P(p_{a,P}) = 0$ for each $k \in I$,
\item $\varphi_{a,P} \pi_{P} = \sum_{k}\nu_{k}g_{a,P}$ for each
$n = 1, \dots, n$ and each $P \in \cat{A}$.
\end{enumerate}
It follows from the third condition that, for every $a$ and $P$,
\begin{equation}\label{una}
\varphi_{a,P} = \varphi_{a,P}\pi_P\iota_P = \sum_{k}
\nu_{k}g_{a,P,k}\iota_P
\end{equation}
Since each $P$ is finitely generated as a right $A$--module, it follows that $g_{a,P,k} : P
\rightarrow \Sigma$ factorizes throughout a finite direct sum $\oplus_{Q \in \cat{F}}Q$. We
can assume that the finite set $\cat{F}$ is independent of $a$ and $P$ (recall that only
finitely many $g_{a,P,k}$'s are nonzero). Therefore,
\begin{equation}
g_{a,P,k}\iota_P = (\sum_{Q \in
\cat{F}}\iota_Q\pi_Q)g_{a,P,k}\iota_P = \sum_{Q \in \cat{F}}
\iota_{Q}\pi_{Q}g_{a,P,k}\iota_P
\end{equation}
In view of \eqref{una}, we have
\begin{equation}\label{fis}
\varphi_{a,P} = \sum_{k,Q} \nu_{k,Q}t_{Q,a,P,k},
\end{equation}
where $\nu_{k,Q} = \nu_{k}\iota_Q \in Q^*$ and $t_{Q,a,P,k} =
\pi_Q g_{a,P,k}\iota_P \in T_{PQ}$.
 On the other hand, $2$ implies that for each $Q,k$
\begin{equation}\label{tes}
\sum_{a,P}t_{Q,a,i,k}p_{a,P} =
\sum_{a,P}\pi_{Q}g_{a,P,k}\iota_{a,P}(p_{a,P}) =
\pi_{Q}(\sum_{a,P}g_{a,P,k}\iota_{a,P}(p_{a,P})) = 0
\end{equation}
Finally,
\[
\begin{array}{lcl}
\sum_{a,P}\varphi_{a,P} \tensor{T_P} p_{a,P} & = &
\sum_{a,P}(\sum_{k,Q} \nu_{k,Q}t_{Q,a,P,k}) \tensor{T_P} p_{a,P}\\
& = & \sum_{a,P,k,Q} \nu_{k,Q}t_{Q,a,P,k} \tensor{T_P}  p_{a,P} -
\sum_{k,Q} \nu_{k,Q} \tensor{T_Q} (\sum_{a,P} t_{Q,a,P,k}p_{a,P})
\\
& = & \sum_{a,k}\left(\sum_{P,Q} \nu_{P,Q}\nu_{k,Q}t_{Q,a,P,k}
\tensor{T_P} p_{a,P} - \sum_{P,Q} \nu_{k,Q}\tensor{T_Q}
t_{Q,a,P,k} p_{a,P} \right),
\end{array}
\]
where the first equality follows from \eqref{fis}, and the second
from \eqref{tes}.
\end{proof}

\begin{definition}
The $A$--coring $\rcomatrix{T}{\Sigma}$ will be called the
\emph{infinite comatrix coring} associated to the category
$\cat{A}$ and the functor $\omega: \cat{A} \rightarrow \add{A_A}$.
Its comultiplication $\Delta$ is given explicitly as follows: once
a (finite) dual basis $\{ e_{\alpha_P}^*, e_{\alpha_P} \}$ is
chosen for each $P \in \cat{A}$ (recall that $P_A$ is finitely
generated and projective), we have, from \eqref{sobreyectivo}
\begin{equation}\label{Deltadef}
\Delta (\varphi \tensor{T} p) = \sum_{ P \in \cat{F}}
\sum_{\alpha_P} \varphi \iota_P \pi_P \tensor{T}
\iota_P(e_{\alpha_P}) \tensor{A} e_{\alpha_P}^*\pi_P \tensor{T}
\iota_P\pi_P (p),
\end{equation}
for $\varphi \tensor{T} p \in \rcomatrix{T}{\Sigma}$, where
$\cat{F}$ is any finite set of objects of $\cat{A}$ such that $p =
\sum_{P \in \cat{F}} \iota_P \pi_P (p)$. The counit $\epsilon$ of
$\rcomatrix{T}{\Sigma}$ is simply the evaluation map $\varphi
\tensor{T} p \mapsto \varphi (p)$.
\end{definition}

\begin{remark}
If ${}_BP_A$ is a $B-A$--bimodule, with $P_A \in \add{A_A}$, then the comatrix coring $P^*
\tensor{B} P$ of \cite[Proposition 2.1]{ElKaoutit/Gomez:2003} is just the infinite comatrix
coring associated to the canonical functor from the additive category $B$ to $\add{A_A}$ which
sends the unique object of $B$ onto $P$.
\end{remark}

Now assume that $\cat{A}$ is a (small) subcategory of the category of right comodules
$\rcomod{\coring{C}}$ over an $A$--coring $\coring{C}$, and that the functor $\omega : \cat{A}
\rightarrow \add{A_A}$ is the restriction of the underlying functor $U : \rcomod{\coring{C}}
\rightarrow \rmod{A}$ (that is, we are taking $\cat{C} = \rcomod{\coring{C}}$). Let $\Sigma =
\bigoplus_{P \in \cat{A}} P$ and $T = \End{\coring{C}}{U}$. A straightforward argument shows
that the functor $\hom{\coring{C}}{\Sigma}{-} : \rcomod{\coring{C}} \rightarrow \rmod{T}$ is
right adjoint to $- \tensor{T} \Sigma : \rmod{T} \rightarrow \rcomod{\coring{C}}$ (the right
comodule structure of $N \tensor{T} \Sigma$ is inherited from $\Sigma$ for each right
$T$--module $N$). The counit of this adjunction evaluated at $\coring{C}$ gives a homomorphism
of $A$--bimodules
\begin{equation}\label{counitatC}
\xymatrix{ \hom{\coring{C}}{\Sigma}{\coring{C}} \tensor{T} \Sigma \ar[r] & \coring{C} & & (f
\tensor{T} u \ar@{|->}[r] & f(u)),}
\end{equation}
which, in conjunction with the canonical isomorphism $\Sigma^* \cong
\hom{\coring{C}}{\Sigma}{\coring{C}} $ (see, e.g. \cite[18.10]{Brzezinski/Wisbauer:2003})
gives a homomorphism of $A$--bimodules
\begin{equation}\label{can}
\xymatrix{ \mathsf{can} :\rcomatrix{T}{\Sigma} \cong \hom{\coring{C}}{\Sigma}{\coring{C}}
\tensor{T} \Sigma \ar[r] & \coring{C}}
\end{equation}

\begin{lemma}\label{candefined}
The map $\mathsf{can} : \rcomatrix{T}{\Sigma} \rightarrow \coring{C}$, defined explicitly by
$\mathsf{can} (\varphi \tensor{T} u) = (\varphi \tensor{A} \coring{C}) \rho_{\Sigma}(u)$ is a
homomorphism of $A$--corings.
\end{lemma}
\begin{proof}
By Proposition \ref{presentacion}, there is a surjective homomorphism of $A$--corings $\Gamma:
\bigoplus_{P \in \cat{A}} \rcomatrix{T_P}{P} \rightarrow \rcomatrix{T}{\Sigma}$. Clearly, it
suffices to prove that $\mathsf{can} \circ \Gamma$ is a homomorphism of $A$--corings, and,
ultimately, that its restriction $\mathsf{can_P}$ to $\rcomatrix{T_P}{P}$ is a homomorphism of
$A$--corings for every $P$. This canonical map $\mathsf{can}_P : \rcomatrix{T_P}{P}
\rightarrow \coring{C}$ coincides with the homomorphism of $A$--corings given in
\cite[Proposition 3]{ElKaoutit/Gomez:2003}.
\end{proof}

We can now apply Gabriel-Popescu's Theorem and state our
Reconstruction Theorem for corings in the following terms.

\begin{theorem}(Reconstruction)\label{reconstruccion}
Let $\coring{C}$ be an $A$--coring and assume that the category $\rcomod{\coring{C}}$ is
abelian, and that there is a generating set $\cat{A}$ of right $\coring{C}$--comodules such
that $P_A$ is finitely generated and projective for every $P \in \cat{A}$. If $\Sigma =
\bigoplus_{P \in \cat{A}}P$ and $T = \End{\coring{C}}{\Sigma}$, then $\mathsf{can} :
\rcomatrix{T}{\Sigma} \rightarrow \coring{C}$ is an isomorphism of $A$--corings.
\end{theorem}
\begin{proof}
If $\rcomod{\coring{C}}$ is abelian, then it is a Grothendieck category \cite[Proposition
1.2]{ElKaoutit/Gomez/Lobillo:2001unp}. Clearly, $\Sigma$ is then a generator for the category
$\rcomod{\coring{C}}$. By Gabriel-Popescu's Theorem \cite{Popesco/Gabriel:1964}, the canonical
map \eqref{counitatC} is an isomorphism. We thus get that $\mathsf{can}$ is an isomorphism of
$A$--corings.
\end{proof}

\begin{remark} The terminology Galois comodule has been introduced
in \cite{Brzezinski/Wisbauer:2003, Brzezinski:unp2003} to refer to a right
$\coring{C}$--comodule $P$ such that $P_A$ is finitely generated and projective and
$\mathsf{can}: \rcomatrix{T}{P} \rightarrow \coring{C}$ is an isomorphism. These comodules
played a role in the characterization of corings having a finitely generated projective
generator \cite[Theorem 3.2]{ElKaoutit/Gomez:2003} (see also the ``Galois comodule structure
theorem'' \cite{Brzezinski:unp2003}), and in the structure theorem for cosemisimple corings
\cite[Theorem 4.4]{ElKaoutit/Gomez:2003}. Theorem \ref{reconstruccion} suggests that it makes
sense to consider \emph{Galois comodules} without finiteness conditions (the important point
here, we think, is that $\rcomatrix{T}{\Sigma}$ is endowed with a coring structure before to
assume that $\mathsf{can}$ is an isomorphism). Another possibility would be to say that
$\cat{A}$ is a \emph{Galois subcategory} of $\rcomod{\coring{C}}$ whenever the coring
homomorphism $\mathsf{can} : \rcomatrix{T}{\Sigma} \rightarrow \coring{C}$ is an isomorphism.
Of course, $\Sigma$ denotes $\bigoplus_{P \in \cat{A}}P$.
\end{remark}

\begin{corollary}\label{caniso}
Let $\coring{C}$ be a coring over a semi-simple artinian ring $A$. Let $\cat{A}$ be a
generating set of finitely generated right $\coring{C}$--comodules. If $\Sigma = \bigoplus_{P
\in \cat{A}}P$ and $T = \End{\coring{C}}{\Sigma}$, then $\mathsf{can} : \rcomatrix{T}{\Sigma}
\rightarrow \coring{C}$ is an isomorphism of $A$--corings.
\end{corollary}
\begin{proof}
Since in this case ${}_A\coring{C}$ is obviously flat, we get that
$\rcomod{\coring{C}}$ is abelian. Moreover, every finitely
generated comodule is finitely generated as a right $A$--module.
Therefore, our coring is in the hypotheses of Theorem
\ref{reconstruccion}.
\end{proof}

\begin{example}
Every coalgebra $C$ over a field $K$ is isomorphic to $\Sigma^*
\tensor{T} \Sigma$, where $\Sigma$ denotes the coproduct of a
generating set of finite-dimensional right $C$--comodules, and $T
= \End{C}{\Sigma}$ (here, $(-)^*$ denotes the $K$--dual functor).
In particular, if $B$ is an algebra over $K$ and $\Sigma$ is the
coproduct of a set of representatives of all finite-dimensional
right $B$--modules, then the finite dual coalgebra $B^0$,
consisting of those $\varphi \in B^*$ such that $Ker \varphi$
contains an ideal $I$ such that $B/I$ is of finite dimension over
$k$, is isomorphic to the coalgebra $\Sigma^* \tensor{T} \Sigma$,
where $T = \End{B}{\Sigma}$.
\end{example}

\section{Corings with a generating set of small
projectives}\label{Freyd}

Let us return to our functor $\omega : \cat{A} \rightarrow
\add{A_A}$. Assume that $\cat{A}$ embeds in an additive category
$\cat{C}$, and that there exists the coproduct $\Sigma =
\bigoplus_{P \in \cat{A}} P$ in the category $\cat{C}$. Assume
further that the functor $\omega : \cat{A} \rightarrow \add{A_A}$
is the restriction of a functor $U : \cat{C} \rightarrow \rmod{A}$
which commutes with the coproduct $\bigoplus_{P \in \cat{A}}P$.
Recall from Section \ref{inf-matrix} that such a category and
functor can be always constructed.

Let $R$ be the twosided ideal of $T = \End{\cat{C}}{\Sigma}$ given by $R = \bigoplus_{P,Q \in
\cat{A}} T_{PQ}$, where $T_{PQ} = \hom{\cat{A}}{P}{Q}$. We will consider $R$ as a ring with a
complete set of pairwise ortogonal idempotents $\{ 1_P : P \in \cat{A} \}$, thought that $R$
has not in general an unit ($1_P$ is the element of $R$ which is the identity of $T_P =
\End{\cat{A}}{P}$ at $P$ and zero elsewhere). Then $\Sigma$ is an $R-A$--bimodule in a
canonical way, and ${\Sigma}^{\dagger} = \oplus_{P \in \cat{A}}P^*$ becomes an
$A-R$--bimodule. We thus get an $A$--bimodule ${\Sigma}^{\dagger} \tensor{R} \Sigma$.

\begin{lemma}
There is a commutative diagram of surjective homomorphisms of
$A$--bimodules
\begin{equation}\label{triangulo}
\xymatrix{\bigoplus_{P \in \cat{A}}\rcomatrix{T_P}{P} \ar^-{\Gamma}[d] \ar^-{\Gamma_1}[rr] & &
{\Sigma}^{\dagger} \tensor{R} \Sigma
\ar^-{\Gamma_2}[dll] \\
\rcomatrix{T}{\Sigma} }
\end{equation}
\end{lemma}
\begin{proof}
For each $P \in \cat{A}$, let $\iota_{P^*} : P^* \rightarrow
{\Sigma}^{\dagger}$ and $\iota_P : P \rightarrow \Sigma$ denote
the canonical inclusions. By $\pi_{P^*} : {\Sigma}^{\dagger}
\rightarrow P^*$ and $\pi_P : \Sigma \rightarrow P$ we denote the
canonical projections. Some straightforward computations show that
the map
\[
\gamma_P : \rcomatrix{T_P}{P} \rightarrow {\Sigma}^{\dagger} \tensor{R} \Sigma \qquad (\varphi
\tensor{T_P} p \mapsto \iota_{P^*}(\varphi) \tensor{R} \iota_P(p))
\]
is a well defined homomorphism of $A$--bimodules. The family $\{ \gamma_P : P \in \cat{A} \}$
determines the homomorphism of $A$--bimodules $\Gamma_1$ in a unique way. In order to show
that $\Gamma_1$ is surjective, observe that every element of ${\Sigma}^{\dagger} \tensor{R}
\Sigma$ is a sum of elements of the form $\iota_{P^*}(\varphi) \tensor{R} \iota_Q(q)$, for
some $\varphi \in P^*, q \in Q$. But if $Q \neq P$, then $\iota_{P^*}(\varphi) \tensor{R}
\iota_Q(q) = \iota_{P^*}(\varphi)\iota_P\pi_P \tensor{R}\iota_Q(q) = 0$, which
proves that $\Gamma_1$ is onto. \\
Now, let us consider the map
\[
\Gamma_2 : {\Sigma}^{\dagger} \tensor{R} \Sigma \rightarrow \rcomatrix{T}{\Sigma} \qquad
\left(\iota_{P^*}(\varphi) \tensor{R} \iota_Q(q) \mapsto \varphi \pi_P \tensor{T} \iota_Q(q),
\; \varphi \in P^*, q \in Q \right)
\]
To check that $\Gamma_2$ is well defined, let $t \in T_{MN} =
\hom{\cat{A}}{M}{N}$, and compute
\[
\begin{array}{lcl}
\varphi \pi_P \iota_N t \pi_M \tensor{T} \iota_Q(p) & = & \varphi
\pi_P \iota_N t \pi_P
\tensor{T} \iota_Q\pi_Q\iota_Q(q) \\
 & = & \varphi  \pi_P \iota_N t \pi_M \iota_Q \pi_Q \tensor{T} \iota_Q(q) \\
 & = & \left\{ \begin{array}{lcl}
               \varphi t \pi_P \tensor{T} \iota_P (q) & \text{if} &
               P=Q=M=N \\
               0 &  & \text{otherwise}
               \end{array} \right.
\end{array}
\]
\[
\begin{array}{lcl}
\varphi  \pi_P \tensor{T}  \iota_N t \pi_M \iota_Q(p) & = &
\varphi \pi_P \iota_P \pi_P
\tensor{T} \iota_N t \pi_M \iota_Q(q) \\
 & = & \varphi  \pi_P \tensor{T} \iota_P \pi_P  \iota_N t \pi_M \iota_Q(q) \\
 & = & \left\{ \begin{array}{lcl}
               \varphi \pi_P \tensor{T} \iota_P t(q) & \text{if} &
               P=Q=M=N \\
               0 &  & \text{otherwise}
               \end{array} \right.
\end{array}
\]
Finally, in the case $P=Q=M=N$, we have
\begin{multline}
\varphi t \pi_P \tensor{T} \iota_P (q) = \varphi \pi_P \iota_P t
\pi_P \tensor{T} \iota_P(q) = \varphi \pi_P \tensor{T} \iota_P t
\pi_P \iota_P (q) = \varphi \pi_P \tensor{T} \iota_P t (q)
\end{multline}
The diagram \eqref{triangulo} is clearly commutative. By Proposition \ref{presentacion},
$\Gamma$ is surjective and, thus, $\Gamma_2$ is so.
\end{proof}

Recall from Lemma \ref{Jota} that the $K$--submodule $\mathfrak{J}$ of $\bigoplus_{P \in
\cat{A}}\rcomatrix{T_P}{P}$ generated by the set $\{ \varphi\tensor{T_Q} tp - \varphi t
\tensor{T_P} p \; : \; \varphi \in Q^*, p \in P, t \in T_{PQ}, P,Q \in \cat{A} \}$ is a
coideal. The factor $A$--coring is denoted by $\mathfrak{R}(\cat{A})$.

\begin{proposition}
The kernel of $\Gamma_1$ is the coideal $\mathfrak{J}$ and, hence, ${\Sigma}^{\dagger}
\tensor{R} \Sigma$ can be endowed with a structure of $A$--coring such that $\Gamma_1$ is an
homomorphism of $A$--corings. In this way, the commutative diagram \eqref{triangulo} induces a
commutative diagram of isomorphisms of $A$--corings
\begin{equation}\label{trianguloiso}
\xymatrix{ \Rcn{\cat{A}} \ar^-{\simeq}[d] \ar^-{\simeq}[rr] & & {\Sigma}^{\dagger} \tensor{R}
\Sigma
\ar^-{\simeq}[dll] \\
\rcomatrix{T}{\Sigma} }
\end{equation}
\end{proposition}
\begin{proof}
Since $\mathfrak{J}$ is already the kernel of $\Gamma$ by Proposition \ref{presentacion} and
$\Gamma_1$ is surjective, it suffices to prove that $\mathfrak{J} \subseteq Ker \Gamma_1$. But
this is a straightforward computation.
\end{proof}

\begin{proposition}\label{Endomorfismo}
Every $P \in \cat{A}$ is a right ${\Sigma}^{\dagger} \tensor{R} \Sigma$--comodule with
structure map
\begin{equation*}
\varrho_P : P \rightarrow P \tensor{A} {\Sigma}^{\dagger} \tensor{R} \Sigma \qquad (p \mapsto
\sum_{\alpha_P} e_{\alpha_P} \tensor{A} \iota_{P^*}(e_{\alpha_P}^*) \tensor{R} \iota_P(p))
\end{equation*}
Consider $S = \bigoplus_{P,Q \in \cat{A}} \hom{A}{P}{Q}$ as a ring (not necessarily with
unit). Then
\begin{equation*}
\hom{{\Sigma}^{\dagger} \tensor{R} \Sigma}{P}{Q} = \{ f \in \hom{A}{P}{Q} \; : \; f \tensor{R}
\iota_P (p) = 1_Q \tensor{R} \iota_Q(f(p)) \; \text{for every } p \in P \}
\end{equation*}
Therefore, the canonical ring extension $R \rightarrow S$ factors through $\bigoplus_{P,Q \in
\cat{A}} \hom{{\Sigma}^{\dagger} \tensor{R} \Sigma}{P}{Q}$.
\end{proposition}
\begin{proof}
A homomorphism of right $A$--modules $f : P \rightarrow Q$ belongs to $\hom{{\Sigma}^{\dagger}
\tensor{R} \Sigma}{P}{Q}$ if and only if
\begin{equation}\label{fcomodulemap}
\sum_{\alpha_Q} e_{\alpha_Q} \tensor{A}
\iota_{Q^*}(e_{\alpha_Q}^*) \tensor{R} \iota_Q(f(p)) =
\sum_{\alpha_P} f(e_{\alpha_P}) \tensor{A} \iota_{P^*} \tensor{R}
\iota_P(p),
\end{equation}
for every $p \in P$. Now,
\begin{equation*}
Q \tensor{A} {\Sigma}^{\dagger} = Q \tensor{A} \bigoplus_{P \in \cat{A}} P^* \cong
\bigoplus_{P \in \cat{A}} P \tensor{A} Q^* \cong \bigoplus_{P \in \cat{A}} \hom{A}{P}{Q},
\end{equation*}
being this last a direct summand, as a right ideal, of $S$. Therefore, $Q \tensor{A}
{\Sigma}^{\dagger} \tensor{R} \Sigma$ is identified with a $K$--submodule of $S \tensor{R}
\Sigma$ for every $Q \in \cat{A}$. With these identifications, equation \eqref{fcomodulemap}
is equivalent to $f \tensor{R} \iota_P (p) = 1_Q \tensor{R} \iota_Q(f(p)) \; \text{for every }
p \in P$.
\end{proof}

Let us now look at the case where $\cat{C} = \rcomod{\coring{C}}$ is the category of right
comodules over an $A$--coring $\coring{C}$, and $\cat{A}$ is a small subcategory of
$\rcomod{\coring{C}}$ whose objects are right $\coring{C}$--comodules which are finitely
generated and projective as right modules over $A$. In this way, $\omega : \cat{A} \rightarrow
\add{A_A}$ is the restriction to $\cat{A}$ of the underlying functor $U : \rcomod{\coring{C}}
\rightarrow \rmod{A}$. The diagram \eqref{trianguloiso} can be now completed to

\begin{equation}\label{diagramacan}
\xymatrix{ \Rcn{\cat{A}} \ar^-{\simeq}[d] \ar^-{\simeq}[rr] & & {\Sigma}^{\dagger} \tensor{R}
\Sigma
\ar^-{\simeq}[dll] \ar_-{\mathsf{can}}[d] \\
\rcomatrix{T}{\Sigma} \ar_-{\mathsf{can}}[rr] & & \coring{C},}
\end{equation}
where the horizontal $\mathsf{can}$ is defined in \eqref{can}, and the vertical one is just
the composite that makes commute the diagram. When $\Sigma$ is a generator for
$\rcomod{\coring{C}}$ both maps are isomorphisms, and the four $A$--corings are isomorphic by
Theorem \ref{reconstruccion}.

The functor $- \tensor{R} T : \rmod{R} \rightarrow \rmod{T}$ has a
right adjoint $\cdot R : \rmod{T} \rightarrow \rmod{R}$ which
sends $X$ onto $XR = \{ xr : r \in R \}$. Consider the diagram of
functors
\begin{equation}\label{diagramafundamental}
\xymatrix{\rmod{A} \ar@<0.5ex>^{- \tensor{A} \coring{C}}[rr] & & \rcomod{\coring{C}}
\ar@<0.5ex>^{U}[ll] \ar@<0.5ex>^{\hom{\coring{C}}{\Sigma}{-}}[rrr]
\ar@<0.5ex>^{\mathcal{F}}[ddrrr] & & & \rmod{T}
\ar@<0.5ex>^{- \tensor{T}\Sigma}[lll] \ar@<0.5ex>^{\cdot R}[dd] \\
& & & & & \\
& & & & & \rmod{R} \ar@<0.5ex>^{- \tensor{R} T}[uu] \ar@<0.5ex>^{ - \tensor{R} \Sigma}[uulll]
}
\end{equation}
where $\mathcal{F} = \hom{\coring{C}}{\Sigma}{-} R$ and $- \tensor{R} \Sigma \simeq -
\tensor{R } T \tensor{T} \Sigma$. Thus, $- \tensor{R} \Sigma$ is left adjoint to
$\mathcal{F}$.

\begin{lemma}\label{natiso}
There are natural isomorphisms $\hom{\coring{C}}{\Sigma}{- \tensor{A} \coring{C}} R \simeq -
\tensor{A} {\Sigma}^{\dagger}$ and $\mathcal{F} \simeq \bigoplus_{P \in \cat{A}}
\hom{\coring{C}}{P}{-}$.
\end{lemma}
\begin{proof}
Since $\hom{\coring{C}}{\Sigma}{- \tensor{A} \coring{C}} \simeq \hom{A}{\Sigma}{-}$ naturally,
we need just to exhibit a natural isomorphism from $ - \tensor{A} {\Sigma}^{\dagger}$ to
$\hom{A}{\Sigma}{-}R$. Now, observe that $- \tensor{A} {\Sigma}^{\dagger} \simeq \oplus_{P \in
\cat{A}}\hom{A}{P}{-}$, and this last functor is easily shown to be naturally isomorphic to
$\hom{A}{\Sigma}{-}R$ via the isomorphism defined at $X \in \rmod{A}$ by the assignment $f
\mapsto f\iota_P\pi_P$ for every $f \in \hom{A}{P}{X}$ and every $P \in \cat{A}$. This last
construction also yields a natural isomorphism $\bigoplus_{P \in \cat{A}}
\hom{\coring{C}}{P}{-} \simeq \mathcal{F}$.
\end{proof}

The left $\coring{C}$--comodule structure of every $P^*$ induces a structure of left
$\coring{C}$--comodule on ${\Sigma}^{\dagger} = \bigoplus_{P \in \cat{A}}P^*$, which will be
used in Proposition \ref{adjuncion2}.

\begin{proposition}\label{adjuncion2}
The functor $ - \tensor{R} \Sigma : \rmod{R} \rightarrow \rmod{A}$ is left adjoint to the
functor $- \tensor{A} {\Sigma}^{\dagger} : \rmod{A} \rightarrow \rmod{R}$, and this adjunction
induces one for right $\coring{C}$--comodules, that is, the functor $- \tensor{R} \Sigma :
\rmod{R} \rightarrow \rcomod{\coring{C}} $ is left adjoint to the functor $-
\cotensor{\coring{C}} {\Sigma}^{\dagger} : \rcomod{\coring{C}} \rightarrow \rmod{R}$. In
particular, $\mathcal{F} \simeq - \cotensor{\coring{C}} {\Sigma}^{\dagger}$.
\end{proposition}
\begin{proof}
The first adjunction follows from Lemma \ref{natiso} and the diagram
\eqref{diagramafundamental}. For the second, let us first observe that $R = R^2$ is a pure
ideal of $T = \End{\coring{C}}{\Sigma}$. Therefore, $R$ is a $T$--coring in a canonical way,
and the category $\rmod{R}$ of unital right $R$--modules is isomorphic to the category
$\rcomod{R}$ of right $R$--comodules. By \cite[Proposition 4.2]{Gomez:2002},
${\Sigma}^{\dagger}$ is a quasi-finite right $R$--comodule and $- \tensor{R} \Sigma$ is left
adjoint to $ - \cotensor{\coring{C}} {\Sigma}^{\dagger}$ (in fact, what we have is that $-
\tensor{R} \Sigma$ is isomorphic to the co-hom functor $h_R({\Sigma}^{\dagger},-)$).
\end{proof}

The coring homomorphism $\mathsf{can} : {\Sigma}^{\dagger} \tensor{R} \Sigma \rightarrow
\coring{C}$ gives a functor $\mathsf{CAN} : \rcomod{{\Sigma}^{\dagger} \tensor{R} \Sigma}
\rightarrow \rcomod{\coring{C}}$.

\begin{proposition}\label{CAN}
Let $\cat{A}$ be a set of right $\coring{C}$--comodules such that every comodule in $\cat{A}$
is finitely generated and projective as a right $A$--module. If $R = \bigoplus_{P,Q \in
\cat{A}} \hom{\coring{C}}{P}{Q}$, then $R = \bigoplus_{P,Q \in \cat{A}}
\hom{{\Sigma}^{\dagger} \tensor{R} \Sigma}{P}{Q}$. Furthermore, we have a commutative diagram
of functors
\begin{equation}\label{CANind}
\xymatrix{\rcomod{\Sigma^{\dagger} \tensor{R} \Sigma} \ar@{->}^-{\mathsf{CAN}}[rr] & &
\rcomod{\coring{C}} \\ & \rmod{R} \ar@{->}_-{-\tensor{R}\Sigma}[ur] \ar@{->}^-{-\tensor{R}
\Sigma}[ul] & }
\end{equation}
\end{proposition}
\begin{proof}
It follows from Lemma \ref{candefined} and \eqref{diagramacan} that right
$\coring{C}$--comodule structure map $\rho_P : P \rightarrow P \tensor{A} \coring{C}$
factorizes as
\[
\xymatrix{P \ar^-{\rho_P}[rr] \ar^-{\varrho_P}[rd] & & P \tensor{A} \coring{C} \\
 & P \tensor{A} {\Sigma}^{\dagger} \tensor{R} \Sigma \ar^-{P \tensor{A} \mathsf{can}}[ru] &}
\]
where the right ${\Sigma}^{\dagger} \tensor{R} \Sigma$--comodule map $\varrho_P$ is defined in
Proposition \ref{Endomorfismo}. This gives, on one hand, that
$\mathsf{CAN}(\Sigma_{{\Sigma}^{\dagger}\tensor{R}\Sigma}) = \Sigma_{\coring{C}}$ and, on the
other hand, that $\hom{{\Sigma}^{\dagger} \tensor{R} \Sigma}{P}{Q} \subseteq
\hom{\coring{C}}{P}{Q}$ for every $P, Q \in \cat{A}$. The converse inclusions
$\hom{\coring{C}}{P}{Q} \subseteq \hom{{\Sigma}^{\dagger} \tensor{R} \Sigma}{P}{Q}$ follow
from Proposition \ref{Endomorfismo}. Finally, $\mathsf{CAN}(Y \tensor{R}
\Sigma_{{\Sigma}^{\dagger} \tensor{R} \Sigma}) = Y \tensor{R} \Sigma_{\coring{C}}$ for every
$Y \in \rmod{R}$, whence the commutativity of the diagram \eqref{CANind}.
\end{proof}

A set $\cat{A}$ of objects of a Grothendieck category $\cat{C}$ is said to be a
\emph{generating set of small projectives} if every object in $\cat{A}$ is small and
$\bigoplus_{P \in \cat{A}}P$ is a projective generator for $\cat{C}$. The following theorem
generalizes \cite[Theorem 3.7]{Schneider:1990}, \cite[Theorem 5.6]{Brzezinski:2002}, and
\cite[Theorem 3.2]{ElKaoutit/Gomez:2003}.

\begin{theorem}\label{descent}
Let $\coring{C}$ be an $A$--coring and $\cat{A}$ a set of right $\coring{C}$--comodules. Let
$\Sigma = \bigoplus_{P \in \cat{A}} P$ and ${\Sigma}^{\dagger} = \bigoplus_{P \in \cat{A}}
P^*$. Consider the ring extension $$R = \bigoplus_{P,Q \in \cat{A}} \hom{\coring{C}}{P}{Q}
\subseteq \bigoplus_{P,Q \in \cat{A}}\hom{A}{P}{Q} = S.$$ The following statements are
equivalent
\begin{enumerate}[(i)]
\item ${}_A\coring{C}$ is flat and $\cat{A}$ is a generating set of small
projectives for $\rcomod{\coring{C}}$;
\item ${}_A\coring{C}$ is flat, every comodule in $\cat{A}$ is
finitely generated and projective as a right $A$--module, $\mathsf{can} : {\Sigma}^{\dagger}
\tensor{R} \Sigma \rightarrow \coring{C}$ is an isomorphism, and $- \tensor{R} \Sigma :
\rmod{R} \rightarrow \rcomod{{\Sigma}^{\dagger} \tensor{R} \Sigma}$ is an equivalence of
categories;
\item every comodule in $\cat{A}$ is
finitely generated and projective as a right $A$--module, $\mathsf{can} : {\Sigma}^{\dagger}
\tensor{R} \Sigma \rightarrow \coring{C}$ is an isomorphism, and ${}_R\Sigma$ is faithfully
flat;
\item ${}_A\coring{C}$ is flat, every comodule in $\cat{A}$ is
finitely generated and projective as a right $A$--module, $\mathsf{can} : {\Sigma}^{\dagger}
\tensor{R} \Sigma \rightarrow \coring{C}$ is an isomorphism, and ${}_RS$ is faithfully flat;
\item ${}_A\coring{C}$ is flat and $-\tensor{R} \Sigma : \rmod{R} \rightarrow
\rcomod{\coring{C}}$ is an equivalence of categories.
\end{enumerate}
\end{theorem}
\begin{proof}
$(i) \Leftrightarrow (v)$ Since ${}_A\coring{C}$ is a flat, it follows from \cite[Proposition
1.2]{ElKaoutit/Gomez/Lobillo:2001unp} that $\rcomod{\coring{C}}$ is a Grothendieck category.
By Proposition \ref{adjuncion2} $- \tensor{R} \Sigma$ is left adjoint to $\cat{F}:
\rcomod{\coring{C}} \rightarrow \rmod{R}$.  By Lemma \ref{natiso}, we can apply Freyd
Theorem's \cite[p. 120]{Freyd:1964} in conjunction with Gabriel Theorem's \cite[Proposition
II.2]{Gabriel:1962} (see also \cite{Harada:1973}) to get that $\cat{F}$ is an equivalence of
categories if and only if $\cat{A}$ is a generating set of small projectives for
$\rcomod{\coring{C}}$. Therefore, $-\tensor{R} \Sigma:\rmod{R} \rightarrow
\rcomod{\coring{C}}$ is itself an equivalence of categories if and only if $\cat{A}$ is a
generating
set of small projectives. \\
$(i) \Rightarrow (ii)$  Since ${}_A\coring{C}$ is flat, $\rcomod{\coring{C}}$ is a
Grothendieck category and the forgetful functor $U_A: \rcomod{\coring{C}} \rightarrow
\rmod{A}$ is exact \cite[Proposition 1.2]{ElKaoutit/Gomez/Lobillo:2001unp}. Moreover, it has
an exact right adjoint $-\tensor{A} \coring{C}$. If $\cat{A}$ is a generating set of small
projectives, then every comodule comodule in $\cat{A}$ is small and projective as a right
$A$--module. By \cite[Section 4.11, Lemma 1]{Pareigis:1970} every $P \in \cat{A}$ is then
finitely generated and projective. On the other hand, Corollary \ref{caniso} and
\eqref{diagramacan} imply that $\mathsf{can} : {\Sigma}^{\dagger} \tensor{R} \Sigma
\rightarrow \coring{C}$ is an isomorphism of $A$--corings, and so $\mathsf{CAN}:
\rcomod{\Sigma^{\dagger} \tensor{R} \Sigma} \rightarrow \rcomod{\coring{C}}$ is already an
equivalence of categories. By Lemma \ref{CAN}, we have that $-\tensor{R}\Sigma: \rmod{R}
\rightarrow
\rcomod{\Sigma^{\dagger} \tensor{R} \Sigma}$ is also an equivalence of categories.\\
$(ii) \Rightarrow (iii)$ The functor $- \tensor{R} \Sigma : \rmod{R} \rightarrow
{\Sigma}^{\dagger} \tensor{R} \Sigma$ is obviously faithful and exact. Since $
{\Sigma}^{\dagger} \tensor{R} \Sigma \cong \coring{C}$ is flat as a left $A$--module, we have,
by \cite[Proposition 1.2 ]{ElKaoutit/Gomez/Lobillo:2001unp}, that the forgetful functor $U :
\rcomod{{\Sigma}^{\dagger} \tensor{R} \Sigma} \rightarrow \rmod{A}$ is faithful and exact.
Therefore, the functor $- \tensor{R} \Sigma : \rmod{R} \rightarrow \rmod{A}$ is faithful and
exact, that is, ${}_R
\Sigma$ is a faithfully flat module. \\
$(iii) \Rightarrow (v)$ Since $\mathsf{can}$ is an isomorphism, ${}_A\coring{C}$ is a flat
module. By Proposition \ref{CAN} we can apply Lemma \ref{adjuncion2} to the infinite comatrix
$A$--coring ${\Sigma}^{\dagger} \tensor{R} \Sigma$ to obtain that the cotensor product functor
$-\cotensor{\Sigma^{\dagger} \tensor{R} \Sigma}{\Sigma}^{\dagger}: \rcomod{\Sigma^{\dagger}
\tensor{R} \Sigma} \rightarrow \rmod{R}$ is right adjoint to the functor $-\tensor{R}\Sigma :
\rmod{R} \rightarrow \rcomod{\Sigma^{\dagger} \tensor{R} \Sigma}$. Since ${}_R\Sigma$ is flat
we have, by a straightforward generalization of \cite[Lemma 2.2]{Gomez:2002} to rings with a
complete set of pairwise ortogonal idempotents, the isomorphism $(M\cotensor{\Sigma^{\dagger}
\tensor{R} \Sigma}{\Sigma}^{\dagger})\tensor{R}\Sigma \cong M \cotensor{\Sigma^{\dagger}
\tensor{R} \Sigma}({\Sigma}^{\dagger}\tensor{R} \Sigma) \cong M$, which turns out to be the
inverse of the counity of the adjunction at $M \in \rcomod{\Sigma^{\dagger} \tensor{R}
\Sigma}$. Moreover, if $\eta_X: X \rightarrow (X \tensor{R} \Sigma)\cotensor{\Sigma^{\dagger}
\tensor{R} \Sigma}{\Sigma}^{\dagger}$ is the unity of the adjunction at $X \in \rmod{R}$, then
an inverse to $\eta_X \tensor{R} \Sigma$ is obtained by the isomorphism $((X \tensor{R}
\Sigma) \cotensor{\Sigma^{\dagger} \tensor{R} \Sigma} {\Sigma}^{\dagger}) \tensor{R} \Sigma
\cong (X \tensor{R} \Sigma) \cotensor{\Sigma^{\dagger} \tensor{R} \Sigma} ({\Sigma}^{\dagger}
\tensor{R} \Sigma) \cong X \tensor{R} \Sigma$. Since ${}_{R} \Sigma$ is faithful, we get that
$\eta_X$ is an isomorphism and, hence $- \tensor{R} \Sigma : \rmod{R} \rightarrow
\rcomod{\Sigma^{\dagger} \tensor{R} \Sigma}$ is an equivalence of categories. It follows from
Lemma \ref{CAN} that $- \tensor{R} \Sigma : \rmod{R} \rightarrow \rcomod{\coring{C}}$ is an
equivalence, as $\mathsf{can}: \Sigma^{\dagger} \tensor{R} \Sigma \rightarrow \coring{C}$ is
an isomorphism. \\
The proof of the equivalence $(iii) \Leftrightarrow (iv)$ is that of \cite[Theorem
3.2]{ElKaoutit/Gomez:2003}, taking the $R$--bilinear isomorphism $S \cong \Sigma \tensor{A}
{\Sigma}^{\dagger}$ into account.
\end{proof}

A consequence of Theorem \ref{descent} is a version for our
functor $\omega: \cat{A} \rightarrow \add{A_A}$ of the Faithfully
Flat Descent Theorem for a (noncommutative) ring extension $B
\subseteq A$ (in categorical words, this is the case when
$\cat{A}$ has a single object whose image under $\omega$ is $A$).

So, let $\omega : \cat{A} \rightarrow \add{A_A}$ be a faithful functor, where $\cat{A}$ is a
small category.  Consider the rings $R = \bigoplus_{P,Q \in \cat{A}} \hom{\cat{A}}{P}{Q}$ and
$\overline{R} =\bigoplus_{P,Q \in \cat{A}} \hom{{\Sigma}^{\dagger} \tensor{R} \Sigma}{P}{Q}$,
and the ring homomorphism $\lambda : R \rightarrow \overline{R}$ defined in Proposition
\ref{Endomorfismo}. Finally, put $\Sigma = \bigoplus_{P \in \cat{A}}P$.

\begin{lemma}\label{can-iso}
 Then the homomorphism of $A$--corings
$\mathsf{can} : {\Sigma}^{\dagger} \tensor{\overline{R}} \Sigma \rightarrow {\Sigma}^{\dagger}
\tensor{R} \Sigma$ is an isomorphism.
\end{lemma}
\begin{proof}
A straightforward computation gives that
$\mathsf{can}(\iota_{P^*}(\varphi)
\tensor{\overline{R}}\iota_P(p)) = \iota_{P^*}(\varphi)
\tensor{R}\iota_P(p)$ for every $\varphi \in P^*, p \in P, P \in
\cat{A}$. So, $\mathsf{can}$ is the inverse to the obvious map
induced by the ring homomorphism $\lambda : R \rightarrow
\overline{R}$.
\end{proof}

\begin{theorem}[Faithfully Flat Descent] With the
previous notations, and \[S = \bigoplus_{P,Q \in
\cat{A}}\hom{A}{P}{Q},
\] the following statements are equivalent.
\begin{enumerate}[(i)]
\item ${}_A(\Sigma^{\dagger} \tensor{R} \Sigma)$ is flat, $\cat{A}$ becomes a generating set of
small projectives for $\rcomod{\Sigma^{\dagger} \tensor{R} \Sigma}$, and $\lambda: R
\rightarrow \overline{R}$ is an isomorphism; \item ${}_A(\Sigma^{\dagger} \tensor{R} \Sigma)$
is flat and $-\tensor{R} \Sigma: \rmod{R} \rightarrow \rcomod{\Sigma^{\dagger} \tensor{R}
\Sigma}$ is an equivalence of a categories;
\item ${}_R\Sigma$ is faithfully flat;
\item ${}_A(\Sigma^{\dagger} \tensor{R} \Sigma)$ is flat (or ${}_R\Sigma$ is flat), and ${}_RS$ is
faithfully flat.
\end{enumerate}
\end{theorem}
\begin{proof}
$(i) \Rightarrow (ii)$ This follows from Lemma \ref{can-iso} and
Theorem \ref{descent}.\\
$(ii) \Rightarrow (i) \text{ and } (iii)$ If ${}_A({\Sigma}^{\dagger} \tensor{R} \Sigma)$ is
flat and $ - \tensor{R} \Sigma : \rmod{R} \rightarrow \rcomod{{\Sigma}^{\dagger} \tensor{R}
\Sigma}$ is an equivalence, then ${}_R\Sigma$ is faithfully flat and the functor $ -
\tensor{R} \Sigma$ sends any projective generator of $\rmod{R}$ onto a projective generator of
the Grothendieck category $\rcomod{{\Sigma}^{\dagger} \tensor{R} \Sigma}$. Therefore $\Sigma
\cong R \tensor{R} \Sigma$ is a projective generator for $\rcomod{{\Sigma}^{\dagger}
\tensor{R} \Sigma}$, since $R$ is a projective generator of $\rmod{R}$ \cite{Gabriel:1962}.
This means that $\cat{A}$ is a generating set of small projectives. By Theorem \ref{descent},
$- \tensor{\overline{R}} \Sigma : \rmod{\overline{R}} \rightarrow \rcomod{{\Sigma}^{\dagger}
\tensor{R} \Sigma}$ is an equivalence of categories. That $\lambda$ is an isomorphism of rings
follows now from the commutative diagram of functors
\[
\xymatrix{\rmod{R} \ar^-{- \tensfun{R} \Sigma}[rr] & &
\rcomod{{\Sigma}^{\dagger} \tensor{R} \Sigma} \\
\rmod{\overline{R}} \ar^-{- \tensfun{\overline{R} } \Sigma}[rr] \ar_-{F}[u] & &
\rcomod{{\Sigma}^{\dagger} \tensor{R} \Sigma} \ar_-{\mathsf{CAN}}[u]},
\]
where $F$ is the restriction scalars functor associated to
$\lambda$, which turns out to be an equivalence of categories. \\
$(iii) \Rightarrow (ii)$ It follows from Proposition \ref{Endomorfismo} that $\lambda
\tensor{R} \Sigma : R \tensor{R} \Sigma \rightarrow \overline{R} \tensor{R} \Sigma$ is an
isomorphism, hence $\lambda$ is an isomorphism as ${}_R\Sigma$ is faithfully flat. The
implication follows
from Lemma \ref{can-iso} and Theorem \ref{descent}.\\
$(iii) \Leftrightarrow (iv)$ The proof of the equivalence between
$(iii)$ and $(iv)$ in Theorem \ref{descent} works here.
\end{proof}

\providecommand{\bysame}{\leavevmode\hbox to3em{\hrulefill}\thinspace}
\providecommand{\MR}{\relax\ifhmode\unskip\space\fi MR }
\providecommand{\MRhref}[2]{%
  \href{http://www.ams.org/mathscinet-getitem?mr=#1}{#2}
} \providecommand{\href}[2]{#2}


\begin{thebibliography}{10}

\bibitem{Brzezinski:1999}
T.~Brzezi{\'n}ski, \emph{{On modules associated to coalgebra Galois
  extensions}}, J. Algebra \textbf{215} (1999), 290--317.

\bibitem{Brzezinski:2002}
T.~Brzezi\'nski, \emph{{The structure of corings: induction functors,
  Maschke-type theorem, and Frobenius and Galois-type properties.}}, Algebr.
  Represent. Theory \textbf{5} (2002), 389--410.

\bibitem{Brzezinski:unp2003}
T.~Brzezi{\'n}ski, \emph{Galois comodules}, Preprint arXiv:math.RA/0312159 v1,
  2003.

\bibitem{Brzezinski/Hajac:1999}
T.~Brzezi{\'n}ski and P.~M. Hajac, \emph{Coalgebra extensions and algebra
  coextensions of {G}alois type}, Comm. Algebra \textbf{27} (1999), no.~3,
  1347--1367.

\bibitem{Brzezinski/Wisbauer:2003}
T.~Brzezi{\'n}ski and R.~Wisbauer, \emph{Corings and comodules}, LMS, vol. 309,
  Cambridge University Press, 2003.

\bibitem{Caenepeel:2003unp}
S.~Caenepeel, \emph{Galois corings from the descent theory point of view},
  preprint arXiv:math.RA/0311377, 2003.

\bibitem{Dade:1980}
E.D. Dade, \emph{Group graded rings and modules}, Math. Z. \textbf{174} (1980),
  241--262.

\bibitem{Doi/Takeuchi:1989}
Y.~Doi and M.~Takeuchi, \emph{{Hopf-Galois extensions of algebras, the
  Miyashita-Ulbrich action, and Azumaya algebras.}}, J. Algebra \textbf{121}
  (1989), 488--516.

\bibitem{ElKaoutit/Gomez:2003}
L.~El~Kaoutit and J.~G\'{o}mez-Torrecillas, \emph{{Comatrix corings: Galois
  corings, descent theory, and a structure theorem for cosemisimple corings.}},
  Math. Z. \textbf{244} (2003), 887--906.

\bibitem{ElKaoutit/Gomez/Lobillo:2001unp}
L.~El~Kaoutit, J.~G{\'o}mez-Torrecillas, and F.~J. Lobillo, \emph{Semisimple
  corings}, preprint arXiv:math.RA/0201070; MPS:Pure mathematics/0202009, 2001.

\bibitem{Freyd:1964}
P.~Freyd, \emph{Abelian categories. {A}n introduction to the theory of
  functors}, Harper's Series in Modern Mathematics, Harper \& Row Publishers,
  New York, 1964.

\bibitem{Gabriel:1962}
P.~Gabriel, \emph{Des cat\'{e}gories ab\'{e}liennes}, Bull. Soc. Math. France
  \textbf{90} (1962), 323--448.

\bibitem{Gomez:2002}
J.~G\'{o}mez-Torrecillas, \emph{{Separable functors in corings.}}, Int. J.
  Math. Math. Sci. \textbf{30} (2002), 203--225.

\bibitem{Guzman:1989}
F.~Guzman, \emph{Cointegrations, relative cohomology for comodules, and
  coseparable corings}, J. Algebra \textbf{126} (1989), 211--224.

\bibitem{Harada:1973}
M.~Harada, \emph{Perfect categories. {II}. {H}ereditary categories}, Osaka J.
  Math. \textbf{10} (1973), 343--355.

\bibitem{Joyal/Street:1990}
A.~Joyal and R.~Street, \emph{An introduction to {T}annaka duality and quantum
  groups}, Category theory (Como, 1990), Lecture Notes in Math., vol. 1488,
  Springer, Berlin, 1991, pp.~413--492.

\bibitem{Koppinen:1994}
M.~Koppinen, \emph{Variations on the smash product with applications to
  group-graded rings}, J. Pure Appl. Algebra \textbf{104} (1994), 61--80.

\bibitem{Kreimer/Takeuchi:1981}
H.F. Kreimer and M.~Takeuchi, \emph{{Hopf algebras and Galois extensions of an
  algebra.}}, Indiana Univ. Math. J. \textbf{30} (1981), 675--692.

\bibitem{Pareigis:1970}
B.~Pareigis, \emph{Categories and {F}unctors}, Academic Press, New York, 1970.

\bibitem{Popesco/Gabriel:1964}
N.~Popesco and P.~Gabriel, \emph{Caract\'erisation des cat\'egories
  ab\'eliennes avec g\'en\'erateurs et limites inductives exactes}, C. R. Acad.
  Sci. Paris \textbf{258} (1964), 4188--4190.

\bibitem{Schneider:1990}
H.~J. Schneider, \emph{Principal homogeneous spaces for arbitrary {H}opf
  algebras}, Israel J. Math. \textbf{72} (1990), 167--195, Hopf algebras.

\bibitem{Stenstrom:1975}
B.~{Stenstr\"om}, \emph{Rings of {Quotients}}, Springer, Berlin, 1975.

\bibitem{Sweedler:1975}
M.~Sweedler, \emph{The predual theorem to the {Jacobson-Bourbaki Theorem}},
  Trans. Amer. Math. Soc. \textbf{213} (1975), 391--406.

\bibitem{Ulbrich:1982}
K.-H. Ulbrich, \emph{{Galoiserweiterungen von nicht-kommutativen Ringen.}},
  Commun. Algebra \textbf{10} (1982), 655--672.

\bibitem{Waterhouse:1979}
W.~C. Waterhouse, \emph{{Introduction to affine group schemes.}}, {Graduate
  Texts in Mathematics. 66. New York, Heidelberg, Berlin: Springer-Verlag. XI,
  164 p. DM 39.50; \$ 22.20 }, 1979.

\end{thebibliography}
\end{document}